\documentclass{article}
\UseRawInputEncoding
\pdfoutput=1
\usepackage{amsmath}
\usepackage{amsthm}
\usepackage{amssymb}
\usepackage{color}
\usepackage{mathrsfs}
\usepackage{graphicx}
\usepackage{amsmath}
\numberwithin{equation}{section}
\def\pl{\partial}

\def\*{\raisebox{.5mm}{*}}

\def\div{{\,\rm div\,}}

\def\i{{\,\rm i\,}}

\def\lam{\lambda}

\def\t#1{{\widetilde#1}}

\def\ol{\overline}

\def\R{\mathbb{R}}

\def\PP{{\cal P}}

\def\b{\beta}

\def\Ga{\Gamma}

\def\Del{\Delta}

\def\<{\langle}
\def\>{\rangle}

\def\Om{\Omega}

\def\Im{{\mbox{\rm Im\,}}}

\def\var{\varphi}

\def\det{\mbox{det\,}}

\def\beq{\arraycolsep=1.5pt\begin{eqnarray}}
\def\eeq{\end{eqnarray}}
\def\bma{\left[\begin{array}}
	\def\ema{\end{array}\right]}
\def\bda{\left|\begin{array}}
	\def\eda{\end{array}\right|}
\def\be{\begin{equation}}
\def\ee{\end{equation}}

\newtheorem{thm}{{}\hskip\parindent Theorem}[section]
\newtheorem{lem}{{}\hskip\parindent Lemma}[section]
\newtheorem{pro}{{}\hskip\parindent Proposition}[section]

\newtheorem{dfn}{{}\hskip\parindent Definition}[section]
\newtheorem{rem}{{}\hskip\parindent Remark}[section]
\begin{document}
 \large
 \title{Inverse problem of recovering the time-dependent damping and  nonlinear terms for wave equations}
 \date{}
 \author{Song-Ren  Fu}
\maketitle
\footnote{}
\footnotetext{
Academy of Mathematics and Systems Science,
Chinese Academy of Sciences, Beijing 100190, P. R.
China;

\ \ School of Mathematical Sciences,
University of Chinese Academy of Sciences, Beijing 100049, P. R.
China;

\ \ e-mail: songrenfu@amss.ac.cn.}

{\bf Abstract.} In this paper, we consider the inverse boundary problems of recovering  the time-dependent nonlinearity and  damping term  for a semilinear wave equation on a Riemannian manifold.  The Carleman estimate and the construction of Gaussian beams together with the higher order linearization are respectively used to derive the uniqueness results of recovering the coefficients.\\


\section{Introduction}
Let $(\Om,g)$ be a Riemannian manifold of dimension $n\ge 2$ with smooth boundary $\pl\Om=\Ga.$ Let $\mathcal M=\Om\times (0,T)$ and $\Sigma=\Ga\times (0,T).$
Assume that $(x,t)=(x_1,\cdots,x_n,x_0=t)$ are local coordinates on $\mathcal M.$
The nonlinear wave equation considered in this paper is given by
\be\label{systembf}
\left\{ \begin{array}{l}
u_{tt}-\Del_g u+b(x,t)u_t+f(x,t,u)=0,\ (x,t)\in \mathcal M,\\
u(x,t)=h(x,t),\quad (x,t)\in \Sigma,\\
u(0,x)=u_0(x),\ u_t(0,x)=u_1(x),\ x\in \Om,
\end{array} \right.
\ee
where $f:\ol\Om\times \R^1\times \mathbb C\rightarrow \mathbb C$ is a smooth function. In the local coordinates, $$\Del_gu=\div_gDu=(\det g)^{-\frac12}\sum\limits_{ij=1}^n\pl_{x_j}[(\det g)^{\frac12}g^{ij}\pl_{x_i}u],$$
where $\div_g$ and $D$ are the divergence operator and Levi-Civita connection in the metric $g$, respectively.

The main goal in this paper is to study the inverse problem of recovering the time-dependent coefficient $b(x,t)$ (damping term), and the nonlinear term $f(x,t,z)$ by some suitable boundary measurements. There are lots of literature concerning the inverse problem of recovering time-dependent or time-independent coefficients in PDEs. We will firstly consider the inverse problem of recovering the particular case, where the nonlinear term $f$ is of time-independent. More precisely, we will recovery $f(x,z)$ in (\ref{systembf}) by means of the Carleman estimates.
It is pointed out that, after the methodology created by \cite{BK}, the Carleman estimates are well used in inverse problems of uniquely determining the time-independent coefficients.

Secondly, for the time-dependent case, we will recovery $b(x,t)$ and $f(x,t,z)$ simultaneously by the higher order linearization method together with constructing some Gaussian beams. We mention that \cite{ISA} devoted to an inverse boundary problem for a nonlinear parabolic equation, in which the first order linearization of the DN map was proposed. This approach has been developed and applied in different other contexts.

For simplicity, we write the inverse problem of recovering $f(x,z)$ as inverse problem (I),  recovering $b(x,t)$ and $f(x,t,z)$ as inverse problem (II).\\

Let $m$ be a positive integer. We introduce the following energy space
$$E^m=\mathop\cap\limits_{k=0}^m C^k([0,T];H^{m-k}(\Om))$$ with the norm $$||u||^2_{E^m}=\mathop{\sup}\limits_{0\le t\le T}\sum\limits_{k=0}^m||\pl^k_tu(t)||^2_{H^{m-k}}\quad {\rm for} \quad u\in E^m.$$
The well-posedness of system (\ref{systembf}) is discussed in the appendix of this present paper.

\subsection{Recovery of the nonlinear term $f(x,z)$}
We consider the following equation.
\be\label{systemlinbq}
\left\{ \begin{array}{l}
u_{tt}-\Del_g u+b(x,t)u_t+f(x,u)=0,\ (x,t)\in \Om\times (0,T),\\
u(x,t)=h(x,t),\quad (x,t)\in \Ga\times (0,T),\\
u(x,0)=\mu(x),\ u_t(x,0)=0,\ x\in \Om.
\end{array} \right.
\ee
Let $u=u(f; \mu)$ be a solution to (\ref{systemlinbq}) with respect to $f$ and $\mu(x).$ Define the input-to-output map as
$$\Lambda_T^f(\mu(x))=\pl_{\nu}u(f;\mu)|_{\Ga\times (0,T)},$$
where $\nu$ denotes the unit normal field pointing outside $\Om$ along $\Ga.$
For simplicity, we assume that $\mu(x)\in C^{\infty}(\ol\Om).$ Based on the initial data $\mu(x),$ the boundary data is given by
$$u(x,t)|_{\Sigma}=h(x,t)=\mu(x)+\frac{t^2}{2}[\Del_g\mu(x)-f(x,\mu(x))]+\sum\limits_{k=3}^m\frac{t^k}{k!}\pl_t^ku(x,0),\quad x\in \Ga,\quad m\ge 2,$$
where 
$$\pl^{k+2}_tu(x,0)=\Del_g(\pl^k_tu(x,0))-\pl_t^{k}(bu_t)|_{t=0}-\pl_t^kf(x,u)|_{t=0}, \quad k=0, 1,\cdots.$$
It is easy to see that the  compatibility conditions hold for system (\ref{systemlinbq}) up to order $m.$

{\bf Main assumptions.}

The following are the main assumptions for the inverse problem (I).

(A.1) $f(x,z)$ is analytic on $\mathbb C$ with values in $C^{\infty}(\mathcal M)$. Moreover, we assume that $f(x,0)=0,$ and $f^{(k)}_z(x,0)|_{\Om_0}=f_{0k}(x)|_{\Om_0},$ where $f_{0k}(x)$ is known for each $k=1, 2,\cdots,$ and $\Om_0$ is an arbitrary small neighborhood of $\Ga$ inside $\Om.$

(A.2) There exists a non-negative strictly convex function $\psi: \ol\Om\rightarrow \R,$
of class $C^3$ in the metric $g.$ There exists a positive constant $\rho$ on $\ol\Om$ such that $\psi(x)$ satisfies

(i) $D^2\psi(x)(X,X)\ge 2\rho|X|_g^2,\quad x\in \Om,\quad X\in \Om_x.$

(ii) $\psi(x)$ has no critical point for $x\in \ol\Om.$ Namely, ${\rm inf}_{x\in\Om}|D\psi(x)|>0.$

(A.3) $b(x,t)\in C^{\infty}(\ol{\mathcal M})$ with $b(x,0)=0.$ Let $\delta>0$ be a small constant. Assume that
$|\mu(x)|\ge \mu_0>0$ for $x\in \ol\Om,$ such that
$$||h||_{H^{m+1}(\Sigma)}+||\mu(x)||_{H^{m+1}(\Om)}\le \frac{\delta}{2},\quad m>\frac n2,$$

We give the admissible sets of functions $b(x,t)$ and $f(x,z)$ as follows. Let $M_0$ be a positive constant. We define
\be\label{admib}
\mathcal{U}_1=\{b(x,t)\in C^{\infty}(\ol{\mathcal M}): b(x,0)=0, \ ||b||_{C^{\infty}(\ol{\mathcal M})}\le M_0\},
\ee
\beq\label{admib}
\mathcal{U}_2=&&\{f(x,z)\in C^{\infty}(\ol\Om\times \mathbb C), f\ {\rm satisfies\ assumption\ (A.1)},\nonumber\\
&&||f^{(k)}_z(x,0)||_{L^{\infty}(\Om)}\le M_0,\ k=1,2,\cdots \}.
\eeq
\begin{rem}
{\rm  Assume that $f$ has the following Taylor expansion
\be
f(x,u)=\sum\limits_{k=0}^{\infty}f^{(k)}_z(x,0)\frac{u^k}{k!}.
\ee
Let $f_1,f_2\in \mathcal U_2$ and let $u_j=u_j(x,t;\mu,f_j)$ be the solution to (\ref{systemlinbq}) with respect to $f_j$ for $j=1,2.$ We see that assumption (A.1) implies that $h(x,t;f_1)=h(x,t;f_2),$
where $h(x,t;f_j)$ is the given  boundary data corresponding to the equation $$u_{jtt}-\Del_gu_j+bu_{jt}+f_j(x,u_j)=0,\quad {\rm for}\ j=1,2.$$

Assumption (A.2) are usually used in the Carleman estimates, see for example \cite{Ymoo,CL}. The existence of convex functions depends on the curvature of $(\Om,g).$ Particularly, for the Euclidean case, we can take $\psi(x)=|x-x_0|^2$ with $x_0\in \R^n\backslash{\ol\Om}.$ For general Riemannian manifolds, such $\psi$ exists locally. There are a number of non-trivial examples to give such $\psi,$ see \cite[Chapter 2.3]{YaoBook}.
}
\end{rem}
We are now in a position to state the main theorem of recovering $f(x,z).$
\begin{thm}
Let assumptions (A.1)-(A.3) hold. Assume that $b(x,t)\in \mathcal U_1$ and $f_1,f_2\in \mathcal U_2.$ Let $T>T^*,$ where $T^*$ is given by (\ref{Tstar}). Then
\be
\Lambda_T^{f_1}(\mu(x))=\Lambda_T^{f_2}(\mu(x))\quad {\rm implies}\quad f_1(x,z)=f_2(x,z),\quad (x,z)\in \Om\times \mathbb C.
\ee
\end{thm}
\begin{rem}
{\rm Assume that $\Gamma_0\subset\Ga$ such that
$$\{x\in \Ga:  \<D\psi,\nu\>\ge 0\}\subset \Ga_0.$$
Then the measurement can be replaced by $\pl_{\nu}u|_{\Ga_0\times (0,T)}.$
}
\end{rem}

\subsection{Recovery of the time-dependent coefficient $b(x,t)$ and the nonlinear term $f(x,t,z)$}
We here consider the nonlinear equation (\ref{systembf}). Let $u(x,t;h,b,f)$ be the solution to (1.1) with respect to $h,b,f.$ With the boundary data, we define the Dirichlet-to-Neumann map as
$$\Lambda^{b,f}_T(h(x,t))=\frac{\pl u(x,t;h)}{\pl\nu}\Big|_{\Sigma}.$$

Before giving the main assumptions for inverse problem (II), we give two definitions as follows.
\begin{dfn}
A Riemannian manifold $(\Om,g)$ is called a simple manifold if it is simply connected, any geodesic in $\Om$ has no conjugate points, and the boundary $\Ga$ is strictly convex with respect to the metric $g$ (the second fundamental form is positive for every point on the boundary).
\end{dfn}
\begin{dfn}
A compact Riemannian manifold $(\Om,g)$ satisfies the foliation condition if there is a smooth strictly convex function, and the boundary $\Ga$ is strictly convex with respect to the metric $g$.
\end{dfn}

The above conditions are crucial in the inverse problems. They are always used as a sufficient condition in the geodesic ray transform. The inversion of the geodesic ray transform on $\Om$ with depends on the geometric properties of $\Om$. See more details in section 4.
Similar to \cite{AFYK}, let $$\mathcal D=\{(x,t)\in \mathcal M: {\rm dist}(x,\Ga)<t<T-{\rm dist}(x,\Ga)\}$$
be the domain of influence. It is clear that no information can be obtained about the coefficients on $\mathcal M\backslash \mathcal D,$ due to the finite speed of propagation of waves.
Let $T>2{\rm Diam(\Om)}$ be given, where $${\rm Diam(\Om)}=\sup \{lengths\ of\ all\ geodesics\ in\ (\Om,g)\}<\infty.$$
We define a subset of $\mathcal D$ by
$$\mathcal E=\{(x,t)\in \mathcal M: D_g(x)<t<T-D_g(x)\},$$
where $D_g(x)$ denotes the length of the longest geodesic passing through the point $x\in \Om.$\\

The following assumptions are the main assumptions for the inverse problem (II).

(B.1) $f(x,t,z)$ is analytic on $\mathbb C$ with values in $C^{\infty}(\mathcal M),$ and $f(x,t,0)=0.$ We define an admissible set for $f$ as
\beq
\mathcal U_3=&&\{f\in C^{\infty}(\ol\Om\times \R\times \mathbb C): f\ {\rm satisfies\ assumption\ (B.1)},\nonumber\\
&& {\rm supp} f^{(k)}_z(x,t,0)\subset\mathcal E,\quad {\rm for}\ k=1,2\}.
\eeq

(B.2) $b(x,t)\in C^{\infty}(\ol{\mathcal M}).$ Similar to $f,$ let $b_0(x,t)$ be a known smooth function.  We define an admissible set for $b$ as
\be
\mathcal U_4=\{b(x,t)\in C^{\infty}(\ol{\mathcal M}):\  {\rm supp}\ b(x,t)\subset\mathcal E\}.
\ee

(B.3) Assume that the initial data $u(x,0)=u_0(x),u_t(x,0)=u_1(x),$ and $$(u_0(x),u_1(x),h(x,t))\in H^{m+1}(\Om)\times H^m(\Om)\times H^{m+1}(\Sigma)$$ satisfying the compatibility conditions (\ref{comcondition}) up to order $m>\frac n2,$  such that
$$||u_0||_{H^{m+1}(\Om)}+||u_1||_{H^m(\Om)}+||h(x,t)||_{H^{m+1}(\Sigma)}\le \frac{\delta}{2},$$
where $\delta>0$ is a given small constant.

Suppose that $(\Om,g)$ is either simple or satisfies the foliation condition, based on the Gaussian beams, which can allow the existence of conjugate points, we have
\begin{thm}
Assume that $u_j=u(x,t;b_j,f_j)$ solves (\ref{systembf}) with respect to $b_j$ and $f_j$ for each $j=1,2.$ Let assumptions (B.1), (B.2) and (B.3) hold, and let $T>2{\rm Diam}(\Om).$ Suppose that $b_1,b_2\in \mathcal U_4,$ $f_1,f_2\in \mathcal U_3$ and $f_z^{(2)}$ is known. Then
$$\Lambda^{b_1,f_1}_T(h(x,t))=\Lambda^{b_2,f_2}_T(h(x,t))\quad {\rm implies}\quad b_1=b_2\ {\rm in}\ \mathcal E,\ {\rm and}\ f_1=f_2\ {\rm in}\ \mathcal D\times\mathbb C.$$
\end{thm}
\subsection{Literature review}
Inverse problems of PDEs have attracted much attention with lots of literature on inverse elliptic equations, parabolic equations, hyperbolic equations, and plate equations, etc. For example, for the inverse elliptic problems, the well known Calder\'{o}n type inverse problems were investigated in \cite{U3,VI,U2,U1} and many subsequent papers. The linear inverse problems of  determining either the time-independent coefficients or the time-dependent coefficients were widely studied. It is known to us that the Carleman estimates method created by \cite{BK} have well used in inverse problems to derive the strong Lipschit stability results for the time-independent coefficients. In this present paper, the Carleman estimate will be used for the recovery of the time-independent nonlinearity $f(x,z)$. However such method deals with the time-independent case only. Also, there are some other methods for the recovery of the time-independent case (or for the real analytic time-dependent coefficients), such as the boundary control (BC) method  steaming from \cite{MB1,MB2} together with Tataru's sharp unique continuation theorem \cite{Tara}. We are not going to list more literature concerning the linear inverse problems.

The inverse problems of nonlinear PDEs are much less. Among them, the early works \cite{ISA,ISA1,ISA2} used the linearization procedure to study the recovery of nonlinear terms appearing in elliptic or parabolic equations. It turns out that the nonlinear interaction of waves can generate new waves, which are essential for the nonlinear inverse problems. For instance, \cite{ALO,OYM, ISA,  HK, ML, Sun} devoted to the unique recovery of nonlinear terms or coefficients appearing in nonlinear elliptic equations. \cite{MC} dealt with the stable recovery of a semilinear term appearing in a parabolic equation, and \cite{RY} studied the fractional semilinear Schr\"odinger equations.

For the inverse problems of hyperbolic equations, we refer to \cite{GN,GN1}, which respectively concerned the the recovery of a conductivity and quadratic coefficients. 
We mention that \cite{MLUW} devoted to the recovery of nonlinear terms from source-to-solution map, where method of the nonlinear interactions of distorted plane waves, originated from \cite{KMU}, was used.
Such method has been successfully used in inverse problems of nonlinear hyperbolic equations, see for example \cite{MHGU,KMU,GUYW,YWTZ}. Similarly, for some semilinear wave equations, instead of the distorted plane waves,  Gaussian beams together with the higher order linearization and the stationary phase method are used to recover the coefficients. For this method, for instance,  we refer to \cite{AFLO,Zhaijian,ZU}. Among those, \cite{Zhaijian} studied an inverse boundary value problem for a semilinear wave equation on a time-dependent Lorentzian manifolds. Both the distorted plane waves and the Gaussian beams were used to derive the uniqueness. In their paper, due to the unsolved inverse problems of recovering the zeroth order term on general manifolds, they assumed the nonlinearity has the following form
$$H(x,z)=\sum\limits_{k=2}^{\infty}h_k(z)z^k,\quad h_k\in C^{\infty}.$$
In their paper, the recovery of the quadric term $f_{uu}$ is much complex and interesting. The distorted plane waves were used to construct four future light-like vectors to recovery the quadric term, but without recovering the zeroth term $f_u(x,t,0).$ However,  the quadric term is unable to recovery by using the Gaussian beams.

It is worth noting that \cite{Kian} considered the inverse problem of determining a general nonlinear term appearing in a semilinear
hyperbolic equation on a Riemannian manifold $(M,g)$ with boundary of dimension $n= 2,3.$  They determined the expression
$F(t,x,u)$ both on the boundary $x\in M$ and inside the manifold $x\in M$ from some partial knowledge of the solutions $u$ on the boundary of the time-space
cylindrical manifold $(0, T)\times M$ or on the lateral boundary $(0,T)\times\pl M.$
Let us also point out that \cite{LLL} investigated inverse boundary problems associated with a time-dependent semilinear hyperbolic equation with variable coefficients. They developed a new method, which combined the observability inequality and a Runge approximation with higher order linearization, to derive the uniqueness of both the initial velocity and the nonlinearity. The measurements they used were either active or passive. For the stability of recovering some coefficients, \cite{MatLas} devoted to the uniqueness and stability of an inverse problem for
a semi-linear wave equation: $u_{tt}-\Del u+a(x,t)u^m=0,$ where $(x,t)\in \R^n\times \R$ and $m\ge 2$. They used the higher order linearization together with the Radon transform to prove the stability results of recovering $a(x,t)$ by the Dirichlet-to-Neumann map. It is worth mentioning that, instead of differentiating the nonlinear equation by $\frac{\pl\epsilon^m}{\pl\epsilon_1\cdots\pl\epsilon_m},$ they used the finite differences operator $D^m_{\epsilon}.$
We note that, all the above literature concerning the nonlinear wave equations do not consider the recovery of the first order coefficients, such as the damping term.  Motivated by these previous works, we study the recovery of the time-dependent damping term and the nonlinearity simultaneously, which can be seen as an extended studying of the existing literature.

The rest of this present paper is organized as follows: In Section 2, we prove Theorem 1.1 by the Carleman setimates. We also add some additional contents in this section.  Section 4 focuses on the Gaussian beams and proofs of Theorem 1.4. Finally, in the appendix, the well-posedness of the semilinear wave equation (\ref{systembf}) is discussed.

\section{Proofs of Theorem 1.1 and an additional result}
In this section, we focus on proving Theorem 1.1 and give an additional result of recovering a leading coefficient.
\subsection{Proofs of Theorem 1.1}
Let $N\ge 1$ be a positive inter, and $\epsilon=(\epsilon_1,\cdots,\epsilon_N)$. Let $\mu_k(x)\ne 0$ for $x\in \ol\Om,$ and
$$\mu(x)=\sum\limits_{k=1}^n\epsilon_k\mu_k(x)\in H^{m+1}(\ol\Om)$$
with $|\epsilon|=\sum\limits_{k=1}^N|\epsilon_k|$ sufficiently small, such that
$$||h(x,t)||_{H^{m+1}(\Sigma)}+||\epsilon_1\mu_1+\cdots+\epsilon_N\mu_N||_{H^{m+1}(\Om)}\le \frac\delta2.$$
By the same output, we have
\be
\Lambda^{f_1}(\epsilon_1\mu_1+\cdots\epsilon_N\mu_N)=\Lambda^{f_2}(\epsilon_1\mu_1+\cdots\epsilon_N\mu_N),\nonumber
\ee
which gives
\be
\frac{\pl^{|\epsilon|}}{\pl\epsilon_1\cdots\pl\epsilon_N}\Big|_{\epsilon=0}\Lambda^{f_1}(\epsilon_1\mu_1+\cdots\epsilon_N\mu_N)
=\frac{\pl^{|\epsilon|}}{\pl\epsilon_1\cdots\pl\epsilon_N}\Big|_{\epsilon=0}\Lambda^{f_2}(\epsilon_1\mu_1+\cdots\epsilon_N\mu_N).
\ee
Clearly, $\Lambda^f(\sum\limits_{k=1}^N\epsilon_k\mu_k)$ contains more information than $\{\Lambda^f(\mu_k)\}_{k=1,\cdots,N}.$ Indeed, useful
information can be extracted from
$$\frac{\pl^N}{\pl\epsilon_1\cdots\pl\epsilon_N}\Big|_{\epsilon=0}\Lambda^f(\sum\limits_{k=1}^N\epsilon_k\mu_k).$$

Let $$\t Q=\Om\times (-T,T),\quad  \t\Sigma=\Gamma\times (-T,T).$$
Extend the domains of $u_j$ and $b(x,t)$ to the region $\t Q$ evenly as usual. Notice that $u_1(x,t)=u_{1t}=b(x,0)=0,$ then $u_{1ttt}(x,0)=0,$ which implies that the extension of $u_1(x,t)$ (also $u_1^j(x,t)$) is smooth.

Similar to \cite[Chapter 5.2]{Ymoo}, we set $$\phi(x,t)=\psi(x)-\beta t^2+\beta_0,\quad \var(x,t)=e^{\lambda\phi(x,t)},\ (x,t)\in \mathcal M,$$
where $\lambda$ is a positive constant, $\beta\in (0,\rho)$ and $\rho$ is given by (i) in assumption (A.2). Moreover,  $\beta_0>0$ is chose such that $\phi(x)>0$. Let
\be\label{Tstar}
T^*=\frac{1}{\sqrt\rho}\left(\mathop{\max}\limits_{x\in \ol\Om}\psi(x)\right)^{\frac12}.
\ee
We assume that $T>T^*.$ Then we can choose $\delta> 0$ and $\beta>0$ such that
$$\rho T^2>\mathop{\max}\limits_{x\in\ol\Om}\psi(x)+4\delta,\quad \beta T^2>\mathop{\max}\limits_{x\in\ol\Om}\psi(x)+4\delta.$$
Thus, $\phi(x,t)$ has the following properties:

(1) $\phi(x,T)\le \beta_0-4\delta$ uniformly for $x\in \Om.$

(2) There exists a small constant $\varepsilon>0$ such that
$$\phi(x,t)\le \beta_0-2\delta\quad {\rm for}\quad (x,t)\in \Om\times [T-2\varepsilon,T]\cup [-T,-T+2\varepsilon].$$
Therefore
$$\var(x,t)\le e^{\lambda(\beta_0-2\delta)}=:d_1<d_0=:e^{\lambda\beta_0}\quad {\rm uniformly\ in}\quad \Om\times [T-2\varepsilon,T].$$

We next give a key Carleman estimate for the linear wave operator $$Lv=v_{tt}-\Del_gv+b(x,t)v_t+q(x,t)v,$$
where $q,b\in L^{\infty}(\mathcal M).$ Let
$$H_0^1(\t Q)=\{u\in H_0^1(\Om\times [-T,T]): u|_{\Ga}=0,\ \pl^l_tu(x,\pm T)=0,\ l=0,1\}.$$
We have the following from \cite{YB}.
\begin{thm}
Under assumptions (A.1) and (A.2), there exist constants $C=C(M_0)>0$ and $\lambda_*>0$, such that for any $\lambda>\lambda_*$, there exists $s_0=s(\lambda)$ such that
\be\label{Car}
\int_{{\t Q}}s[(|Dv|+v_t^2)+s^2v^2]e^{2s\varphi}dgdt\le C\int_{\t Q}|Pv|^2e^{2s\varphi}dgdt+C\int_\Sigma|\pl_{\nu}v|^2e^{2s\var}d\Sigma
\ee
holds for all $v\in H_0^1(\t Q)$ and $s>s_0>1.$
\end{thm}

{\bf Proof of Theorem 1.1.}
Based on the above linearization procedure and the Carleman estimate, we divide the proofs into two steps.

{\bf Step 1.}\ First order linearization. Let $u_j=u_j(x,t;f_j,\epsilon)\in E^{m+1}$ be a solution to (\ref{systemlinbq}) with respect to $f_j$ and $\mu(x)=\epsilon_1\mu_1(x)+\cdots+\epsilon_N\mu_N(x)$ for $j=1,2.$ Then $\t u_j=u_j(x,t; f_j,0)$ solves
\be
\left\{ \begin{array}{l}
\t u_{jtt}-\Del_g \t u_j+b(x,t)\t u_{jt}+f_j(x,\t u_j)=0,\ (x,t)\in \Om\times (0,T),\\
\t u_j(x,t)=0,\quad (x,t)\in \Ga\times (0,T),\\
\t u_j= \t u_{jt}(x,0)=0,\quad x\in \Om,
\end{array} \right.
\ee
which admits a zero solution $\t u_j=0$ since $f_j(x,0)=0$ for each $j=1,2.$
We next linearize the system (\ref{systemlinbq}) around $\t u_j=0$.

Let $N=1$ and let $u_1^j=\frac{\pl}{\pl\epsilon_1}\Big|_{\epsilon=0}u_j.$ Then $u_1^j$ satisfies
\be\label{sysuj}
\left\{ \begin{array}{l}
u^j_{1tt}-\Del_g u^j_1+b(x,t)u^j_{1t}+f_{ju}(x,0)u^j_1=0,\quad (x,t)\in \t Q,\\
u^j_1(x,t)=\frac{\pl}{\pl\epsilon_1}\Big|_{\epsilon=0}h:=h_1,\quad (x,t)\in \t\Sigma,\\
u^j_1(x,0)=\mu_1(x),\ u^j_{1t}(x,0)=0,\quad x\in \Om.
\end{array} \right.
\ee
Set $u_1=u_1^1-u_1^2$ and $q(x)=f_{2u}(x,0)-f_{1u}(x,0),$ then
\be
\left\{ \begin{array}{l}
u_{1tt}-\Del_g u_1+b(x,t)u_{1t}+f_{1u}(x,0)u_1=q(x)u_1^2,\quad (x,t)\in\t Q,\\
u_1(x,t)=0,\quad (x,t)\in \t\Sigma,\\
u_1(x,0)=u_{1t}(x,0)=0,\quad x\in \Om.
\end{array} \right.
\ee
Let $y_1=u_{1t}.$ Then
\be\label{y1}
\left\{ \begin{array}{l}
y_{1tt}-\Del_g y_1+by_{1t}+(b_t+f_{1u}(x,0))y_1=q(x)u_{1t}^2,\quad (x,t)\in \t Q,\\
y_1(x,t)=0,\quad (x,t)\in \t\Sigma,\\
y_1(x,0)=0,\ y_{1t}(x,0)=q(x)\mu_1(x),\quad x\in \Om.
\end{array} \right.
\ee
Next, we chose a cut-off function $\chi(t)\in C_0^{\infty}([-T,T])$ satisfying
\be
0\leq\chi(t)\leq1,\quad
\chi (t) = \left\{ \begin{array}{l}
0,\quad t \in [-T,-T+\varepsilon) \cup (T-\varepsilon ,T],\vspace{1.0ex}\\
1,\quad t \in [-T+2\varepsilon, T-2\varepsilon].
\end{array} \right.\label{4.9}
\ee
Let $\hat y_1=\chi(t)y_1.$ Then
$$\hat y_{1tt}-\Del_g\hat y_1+b\hat y_{1t}+(b_t+f_{1u}(x,0))\hat y_1=\chi q(x)u_{1t}^2+\chi_t(2y_{1t}+by_1)+\chi_{tt}y_1.$$
Since the assumption that $|\mu_1(x)|\ge c_1>0,$ then
\beq\label{esq}
&&\int_{\Om}|q(x)|^2e^{2s\var(x,0)}dg\le C\int_\Om|q(x)\mu_1(x)|^2e^{2s\var(x,0)}dg\nonumber\\
&&=-\int_0^T\frac{\pl}{\pl t}\int_\Om |y_{1t}(x,t)\chi(t)|^2e^{2s\var}dtdg\nonumber\\
&&=-2\int_Q(\chi^2y_{1t}y_{1tt}+\chi\chi_ty_{1t}^2+s\var_t\chi^2y_{1t}^2)e^{2s\var}dgdt\nonumber\\
&&\le Cs\int_Q[|(\chi y_1)_t|^2+\chi_t^2(y_1^2+y_{1t}^2)]e^{2s\var}dgdt+2\int_Q\chi^2y_{1t}y_{1tt}e^{2s\var}dgdt.
\eeq
We compute the term $\chi^2y_{1t}y_{1tt}e^{2s\var}$ as follows.
\beq
\chi^2y_{1t}y_{1tt}e^{2s\var}&&=\chi^2y_{1t}[\Del_g y_1-by_{1t}-(b_t+f_{1u}(x,0))y_1+q(x)u_{1t}^2]e^{2s\var}\nonumber\\
&&=\chi^2\div(y_{1t}e^{2s\var}Dy_1)-\frac12(\chi^2e^{2s\var}|Dy_1|^2)_t+(\chi\chi_t+s\chi^2\var_t)|Dy_1|^2e^{2s\var}\nonumber\\
&&\quad-2s\chi^2y_{1t}\<Dy_1,D\var\>e^{2s\var}-\chi^2[(b_t+f_{1u}(x,0))y_1-q(x)u_{1t}^2]e^{2s\var}\nonumber\\
&&\le \chi^2\div(y_{1t}e^{2s\var}Dy_1)-\frac12(\chi^2e^{2s\var}|Dy_1|^2)_t\nonumber\\
&&\quad+Cs[|D(\chi y_1)|^2+|(\chi y_1)_t|^2+\chi_t^2(y_{1t}^2+|Dy_1|^2)]e^{2s\var}+\chi^2q(x)u_{1t}^2e^{2s\var}.\nonumber
\eeq
Moreover, for system (\ref{sysuj}), the hyperbolic regularity (see, Lemma A.1 in the appendix) implies that $u_{1t}^2\in C([-T,T];H^m(\Om))$ for $m>\frac n2.$ Then the Sobolev embedding theorem shows that $u_{1t}^2\in L^{\infty}(\t Q).$ Thus, there exists a positive constant $C=C(M_0)$ such that
\beq\label{ytytt}
&&-\int_Q\{\chi^2y_{1t}[by_{1t}+(b_t+f_{1u}(x,0))y_1-q(x)u_{1t}^2]\}e^{2s\var}dgdt\nonumber\\
&&\le C\int_Q[\chi_t^2y_1^2+|(\chi y_1)_t|^2+q(x)]e^{2s\var}dgdt.\nonumber
\eeq
Notice that $\frac12(\chi^2|Dy_1|^2e^{2s\var})|_0^T=0.$ Therefore, by (\ref{esq}) and (\ref{ytytt}), with $v=\hat y,$ it follows from the Carleman estimate (\ref{Car}) that
\beq
&&\int_\Om q^2e^{2s\var(x,0)}dg\le Cs\int_Q(\hat y^2+|D\hat y|^2+\hat y_t^2)e^{2s\var}dgdt\nonumber\\
&&\quad+Cs\int_Q\chi_t^2(y_1^2+y_{1t}^2+|Dy_1|^2)e^{2s\var}dgdt+C\int_Q q^2e^{2s\var}dgdt\nonumber\\
&&\le C\int_Qq^2e^{2s\var}dgdt+\int_Q(\chi_t^2+\chi_{tt}^2)(y_1^2+y_{1t}^2+|Dy_1|^2)e^{2s\var}dgdt+Ce^{Cs}||\pl_{\nu}y_1||^2_{L^2(\Sigma)}.\nonumber
\eeq
By the standard energy estimate of system (\ref{y1}), there exists a positive constant $C=C(T,M_0)$ such that
\be\label{energyest}
\int_\Om(y_1^2+y_{1t}^2+|Dy_1|^2)dg\le C\int_\Om q^2dg+C\int_\Sigma|\pl_{\nu}y_1|^2d\Sigma.
\ee
Since $\chi_t,\chi_{tt}\ne 0$ in the case where $\var(x,t)\le d_1,$ with (\ref{energyest}), we have
\be
\int_\Om q^2e^{2s\var(x,0)}dg\le C\int_Q q^2e^{2s\var}dgdt+Ce^{Cs}||\pl_{\nu}y_1||^2_{L^2(\Sigma)}.
\ee
Moreover, by the Lebesgue's theorem, we have
$$\int_Q q^2e^{2s\var}dgdt=\int_\Om q^2e^{2s\var(x,0)}\int_0^Te^{-2s(\var(x,0-\var(x,t)))}dt=o(1)\int_\Om q^2e^{2s\var(x,0)}dg.$$
Thus, $$||q||^2_{L^2(\Om)}\le Ce^{Cs}||\pl_{\nu}y_1||^2_{L^2(\Sigma)}=Ce^{Cs}||\pl_{\nu}(u_1^1-u_1^2)_t||^2_{L^2(\Sigma)}.$$
Hence  $\pl_{\nu}u^1_1=\pl_{\nu}u_1^2$ on $\Sigma$ implies that
$f_{1u}(x,0)=f_{2u}(x,0)=f_u(x,0)$ for $x\in \Om.$\\

{\bf Step 2.} Higher order linearization. Based on Step 1, we set $f_{1u}(x,0)=f_{2u}(x,0)=f_u(x,0)$ for simplicity. Let $u_2^j=\frac{\pl^2}{\pl\epsilon_1\pl\epsilon_2}\Big|_{\epsilon =0}u_j$ for $j=1,2.$ Then
\be
\left\{ \begin{array}{l}
u^j_{2tt}-\Del_g u^j_2+b(x,t)u^j_{2t}+f_u(x,0)u^j_2+f_{juu}(x,0)u_1^1u_1^2=0,\ (x,t)\in \Om\times (0,T),\\
u^j_2(x,t)=0,\quad (x,t)\in \Ga\times (0,T),\\
u^j_2(x,0)=0,\ u^j_{2t}(x,0)=0,\ x\in \Om.
\end{array} \right.
\ee
Here $u_1^1$ and $u_1^2$ satisfy
\be
\left\{ \begin{array}{l}
u^j_{1tt}-\Del_g u^j_1+b(x,t)u^j_{1t}+f_u(x,0)u_1^j=0,\ (x,t)\in \Om\times (0,T),\\
u^j_1(x,t)=\frac{\pl}{\pl\epsilon_j}\Big|_{\epsilon=0}:=h_j(x,t),\quad (x,t)\in \Ga\times (0,T),\\
u^j_1(x,0)=\mu_j(x),\ u^j_{1t}(x,0)=0,\quad x\in \Om.
\end{array} \right.
\ee
Set $y_2=(u_2^1-u_2^2)_t,$  then
\be
\left\{ \begin{array}{l}
y_{2tt}-\Del_g y_2+b(x,t)y_{2t}+(b_t+f_u(x,0))y_2=F(x,t),\ (x,t)\in \Om\times (0,T),\\
y_2(x,t)=0,\quad (x,t)\in \Ga\times (0,T),\\
y_2(x,0)=0,\ y_{2t}(x,0)=-(f_{1uu}-f_{2uu})(x,0)\mu_1(x)\mu_2(x),\quad x\in \Om,
\end{array} \right.
\ee
where $F(x,t)=-(f_{1uu}-f_{2uu})(x,0)(u_{1t}^1u_1^2+u_1^1u_{1t}^2).$
By a similar argument as that in Step 1, we have
$$\pl_{\nu}u_2^1=\pl_{\nu}u_2^2\quad {\rm on}\quad \Sigma\Rightarrow f_{1uu}(x,0)=f_{2uu}(x,0)\Rightarrow u_2^1=u_2^2:=u_2.$$

Suppose that
$$f_{1u}^{(N-1)}(x,0)=f_{2u}^{(N-1)}(x,0),\quad x\in \Om.$$
Let $w^N_j=\frac{\pl^N}{\pl\epsilon_1\cdots\epsilon_N}\Big|_{\epsilon=0}u_j.$
By the recursive assumption, for the expression
$$f_j(x,u)=\sum\limits_{k=0}^{\infty}f^{(k)}_u(x,0)\frac{u^k}{k!},$$ and $v_k=\frac{\pl}{\pl\epsilon_k}\Big|_{\epsilon=0}u,$ we know that
$$\frac{\pl^N}{\pl\epsilon_1\cdots\epsilon_N}\Big|_{\epsilon=0}f_j(x,u)-f_u(x,0)w_j^N-f^{(N)}_{ju}(x,0)v_1\cdots v_N$$
is already known.

The above procedure can be proceeded for $N$-th linearization to obtain
$$f_{1u}^{(N)}(x,0)=f_{2u}^{(N)}(x,0),\quad x\in \Om.$$
Up to now, by the analyticity of $f,$ we have
$$f_1(x,z)=\sum\limits_{k=1}^\infty f^{(k)}_{1u}(x,0)\frac{z^k}{k!}=\sum\limits_{k=1}^\infty f^{(k)}_{2u}(x,0)\frac{z^k}{k!}=f_2(x,z).$$

Thus, the uniqueness result in Theorem 1.1 holds and the proof of Theorem 1.1 is complete. \qedsymbol
\subsection{An additional result of recovering a leading coefficient}
Based on the Carleman estimate, we introduce here a stability result of recovering a leading coefficient.
Suppose that there is a leading coefficient appearing in the wave equation
\be\label{additional}
\varrho(x)u_{tt}-\Del_gu+bu_t+f(x,t,u)=0,\quad (x,t)\in \mathcal M.
\ee
The unique recovery of the mass density $\varrho(x)$ for a wave equation is essential in inverse problems. There are literature concerning this topic, see for example \cite{LBMC,SLLO,PG3}.
We give a brief discussion of recovering $\varrho(x).$

Due to the presence of $\varrho(x),$ we consider a new metric $\hat g=\varrho g,$ then
$$\Del_{\hat g}u=\frac1\varrho\Del_gu+\frac{n-2}{2\varrho^2}\<D\varrho,Du\>.$$
Let $\hat D$ be the Levi-Civita connection in the metric $\hat g.$ It is well known that, for any vectors $X,Y$
$$\hat D_XY=D_XY+\frac12g(D\ln\varrho,Y)X+\frac12g(D\ln\varrho,X)Y+\frac12g(X,Y)D\ln\varrho.$$
If $\psi(x)$ is strictly convex in the metric $\hat g$, then one needs
\beq
&&\hat D^2\psi(X,X)=\varrho g(\hat D_X\hat D\psi,X)=\varrho[g(D_X(\varrho^{-1}D\psi), X)]\nonumber\\
&&\quad+\frac12g(D\psi(\ln\varrho)X+X(\ln\varrho)D\vartheta+X(\psi)D(\ln\varrho),X)\nonumber\\
&&=D^2\psi(X,X)+\frac12D\vartheta(\ln\varrho)|X|_g^2\ge (2+\frac12D\psi(\ln\varrho))|X|_g^2\nonumber\\
&&=\frac1\varrho(2+\frac12D\psi(\ln\varrho))|X|_{\hat g}^2\ge \vartheta_0|X|_{\hat g}^2,\nonumber
\eeq
where $\vartheta_0>0$ is a constant. Therefore, we need further assumption on $\varrho(x)$:
$$g(D\psi,D\varrho)\ge 2\varrho(\vartheta_0-2).$$
Let $u_j(x,t)=u(x,t;\rho_j)$ be a solution to (\ref{additional}) with respect to $\varrho_j$ for each $j=1,2.$ Let $w=u_1-u_2,$ $\hat\varrho=\rho_1-\rho_2.$ Then
\be
\varrho_2 w_{tt}-\Del_gw+bw_t+cw=-\hat\varrho u_{1tt},
\ee
where $c(x,t)=\int_0^1f_u(x,t,ru_1+(1-r)u_2)dr\in L^{\infty}(\mathcal M)$ for sufficiently smooth $u_1$ and $u_2.$ Suppose that $|\Del_g\mu(x)-f(x,0,\mu(x))|>0$ for $x\in \Om,$ where $(u_j(x,0), u_{jt}(x,0))=(\mu(x),0)$ for $j=1,2.$ Then a similar argument with the proof of Theorem 1.1 yields the following stability of recovering $\varrho(x)$:
$$||\varrho_1(x)-\varrho_2(x)||_{L^2(\Om)}\le C||\pl_{\nu}(u_1-u_2)_t||_{L^2(\Sigma)}.$$
\begin{rem}
{\rm
Clearly, as we have discussed above, the non-degeneracy of initial data $u(x,0),$ which is viewed as the input is needed. Moreover, the existence of some strictly convex functions, which seems not sharp, should be assumed. Such functions guarantee that interior information of solutions
to the system arrives at boundary in a finite time.
}
\end{rem}

\section{Gaussian beams and proofs of Theorem 1.4}
In this section, as usual, the boundary data $h$ is the input. A geodesic $\beta(t)\subset\mathcal M$ is called a null geodesic in the metric $\ol g=-dt^2+g$, if $D^{\ol g}_{\dot\beta}\dot\beta=\ol g(\dot\beta,\dot\beta)=0,$ where $D^{\ol g}$ is the connection in the metric $\ol g.$  We intend to construct some Gaussian beams around a null geodesic. For completeness, we give the details of constructing the Gaussian beams.

\subsection{Gaussian beams}
Let $\beta(t)$ be a null geodesic. We firstly introduce the Fermi coordinates in a neighborhood of the null geodesic $\beta$. We follow the constructions in \cite{AJYL}, see also \cite{AFLO}. Recall that $\Om\subset\subset\Om_1$ and the functions are extended smoothly to $\Om_1.$ Let $\beta(t)=(t,\gamma(t))\subset\R\times\Om_1,$ where $\gamma(t)$ is a unit-speed geodesic in the Riemannian manifold $(\Om_1,g).$
Assume that $\beta(t)$ passes through a point $(t_0,x_0),$ where $t_0\in(0,T)$ and $\gamma(t_0)=x_0\in \Om.$ Let $\beta$ join two points $(t_-,\gamma(t_-))$ and $(t_+,\gamma(t_+))$ with $t_{\pm}\in (0,T)$ and $\gamma(t_{\pm})\in \Ga.$
We extend $\beta$ to a larger manifold $\mathcal{M}_1=(0,T)\times\Om_1$ such that $\gamma(t)$ is well defined on $[t_--\epsilon,t_++\epsilon]\subset (0,T)$
with $\epsilon>0$ sufficiently small.

Since the geodesic $\gamma$ is parallel along itself, we can choose $\{e_2,\cdots,e_n\}$ such that $\{\dot\gamma(t_0),e_2\cdots,e_n\}$ forms an orthonormal basis of $\Om_{x_0}.$ Let $s$ denote the arc length along $\gamma$ from $x_0.$ 
Let $E_k(s)\in \Om_{\gamma(s)}$ be the parallel transport of $e_k$ along $\gamma$ to the point $\gamma(s).$

We now define a map
$$\mathcal{F}_1:\R^{n+1}\to \mathcal M_1$$ such that
$$\mathcal F_1(y_0=t;s=y_1,y_2\cdots, y_n)=(t,\exp_{\gamma(s)}(y_2E_2(s)+\cdots+y_nE_n(s))),$$
where $\exp_p(\cdot)$ denotes the exponential map on $\Om_1$ at the point $p$.
In the new coordinates, we have
$$g|_{\gamma}=\sum\limits_{k=1}^ndy_k^2\quad {\rm and}\quad \frac{\pl g_{ij}}{\pl y_k}\Big|_{\gamma}=0\quad {\rm for}\ 1\le i,j,k\le n.$$
On the Lorentzian manifold $(\mathcal{M}_1,-dt^2+g),$ we introduce the well known Fermi coordinates near the null geodesic $\beta: (t_--\epsilon,t_++\epsilon)\rightarrow\mathcal M_1$ as follows.

Let
$$a=\sqrt2(t_--\epsilon),\quad b=\sqrt2(t_++\epsilon),\quad a_0=\sqrt2(t_--\frac{\epsilon}{\sqrt2}),\quad b_0=\sqrt2(t_++\frac{\epsilon}{\sqrt2}),$$ and
$$z_0=\frac{1}{\sqrt2}(t+s)+\frac a2,\quad z_1=\frac{1}{\sqrt2}(-t+s)+\frac a2,\quad z_j=y_j\quad {\rm for}\ 2\le j\le n.$$
Then, we have
\be\label{nullg}
\ol g|_{\beta}=2dz_0dz_1+\sum\limits_{k=2}^ndz_k^2\quad {\rm and}\quad \frac{\pl \ol g_{ij}}{\pl z_k}\Big|_{\beta}=0\quad {\rm for}\ 0\le i,j,k\le n.
\ee
For simplicity, we use the notation $z=(z_0,z^{'})=(z_0,z_1,z^{''})$ to denote the so called Fermi coordinates.  The following lemma from \cite[Lemma 3.1]{AJYL} (see also \cite[Lemma 1]{AFLO}) is essential for the construction of Gaussian beams.
\begin{lem}
Let $\beta:(t_--\epsilon,t_++\epsilon)\rightarrow \mathcal M_1$ be a null geodesic as above. Then there exists a coordinate neighborhood $(U,\Phi)$ of $\beta(t_--\frac\epsilon2,t_++\frac\epsilon2),$ with the coordinates denoted by $(z_0,z^{'})$ such that
$$V=\Phi(U)=(a,b)\times B(0,\delta),$$
where $B(0,\delta)$ denotes a ball in $\R^n$ with a small radius $\delta.$
\end{lem}
Based on the above coordinates, we will construct some approximate Gauss beams in a neighborhood of $\beta$ by defining
$$\mathcal V=\{(z_0,z^{'})\in \mathcal M_1: z_0\in (a_0,b_0), |z^{'}|\le \delta^{'}\}$$ with $0<\delta^{'}<\delta$ sufficiently small such that the set $\mathcal V$ does not
intersect the sets $\{0\}\times \Om$ and $\{T\}\times \Om.$

We use the shorthand notation $\mathcal L_{b,q}u=u_{tt}-\Del_gu+bu_t+qu.$ We consider the WKB ansatz
$$u_\sigma=e^{\i\sigma\varphi}a+r_\sigma,$$
where $\sigma>0$ is a constant, $r_\sigma$ is the reminder term. $\varphi$ and $a$ are called the amplitude and phase respectively. In particular, we will construct $\varphi\in C^{\infty}(\mathcal V)$ and $a\in C_0^{\infty}(\mathcal V).$

Directly calculation yields
\beq\label{lbq}
&&\mathcal L_{b,q}(e^{\i\sigma\var}a)=e^{\i\sigma\var}\mathcal L_{b,q}a+2\i\sigma(\var_ta_t-\<Da,D\var\>_g)e^{\i\sigma\var}\nonumber\\
&&\quad+\i\sigma a(\var_{tt}-\Del_g\var+b\var_t)e^{\i\sigma\var}+a\sigma^2(|D\var|_g^2-\var_t^2)e^{\i\sigma\var}.
\eeq
Based on the above equation, we respectively solve

The eikonal equation $$\mathcal{S}\var=|D\var|_g^2-\var_t^2=\<d\var,d\var\>_{\ol g}=\sum\limits_{k,l=0}^n\ol g^{kl}\pl_k\var\pl_l\var=0,$$
and the
transport equation
\beq\label{transport0}
\mathcal T_b(a,\var)&&=-2\<d\var,da\>_{\ol g}-(\Del_{\ol g}\var-b\var_t)a\nonumber\\
&&=2(\var_ta_t-\<Da,D\var\>_g)+(\var_{tt}-\Del_g\var+b\var_t)a=0.
\eeq
To achieve this, we make the following ansatz for $\var,a,$ namely
$$\var(z_0,z^{'})=\sum\limits_{k=0}^N\var_k(z_0,z^{'}),\quad \quad a=\sum\limits_{k=0}^N\sigma^{-k}\chi(\frac{|z^{'}|}{\delta^{'}})a_k(z_0,z^{'}),\quad a_k=\sum\limits_{j=0}^Na_{k,j}(z_0,z^{'}).$$
Here $\var_k$ is a homogeneous polynomial of degree $k$ with respect
to the variables $z_i$ for $k=0,1,\cdots,N.$
In terms of the Fermi coordinates $z=(z_0,z_1,\cdots,z_n)$ for $\ol g=-dt^2+g,$ we need
\be\label{phase0}
\frac{\pl^{|\Theta|}}{\pl z^\Theta}(\mathcal S\var)(z_0,0)=\frac{\pl^{|\Theta|}}{\pl z^\Theta}\<d\var,d\var\>_{\ol g}|_{\beta}=0\quad {\rm for}\quad \Theta=(0,\theta_1,\cdots,\theta_n),\quad z_0\in (a_0,b_0),\nonumber
\ee
where $\theta_j\ge 0$ are integers for $1\le j\le n,$  and $|\Theta|=\sum\limits_{j=1}^n\theta_j\le N,$ Moreover, for $k=1,\cdots,N,$ we need
$$\frac{\pl^{|\Theta|}}{\pl z^\Theta}\mathcal T(a_0,\var)(z_0,0)=0,\quad \frac{\pl^{|\Theta|}}{\pl z^\Theta}(\i \mathcal T(a_k,\var)+\mathcal L_{b,q}a_{k-1})(z_0,0)=0, \quad z_0\in (a_0,b_0).$$
\vspace{0.1cm}

\noindent{\bf Construction of the phase.}\quad We firstly solve equation (\ref{phase0}) with $|\Theta|=0.$ That is,
$$\sum\limits_{k,l=0}^n\ol g^{kl}\pl_k\var\pl_l\var=0\quad {\rm on}\quad \beta.$$
By (\ref{nullg}), this reduces to
\be\label{var0}
2\pl_0\var\pl_1\var+\sum\limits_{k=2}^n(\pl_k\var)^2=0.
\ee
Similar, for $|\Theta|=1,$ we have
\be\label{var1}
\sum\limits_{kl=0}^n\pl^2_{jk}\var\pl_l\var=0\quad {\rm for}\quad 1\le j\le n.
\ee
Clearly, equations (\ref{var0}) and (\ref{var1}) are satisfied  by respectively setting
$$\var_0=0,\quad \var_1=z_1=\frac{1}{\sqrt2}(-t+s)+\frac a2.$$
For the case where $|\Theta|=2,$ we set
\be\label{var2}
\var_2(z_0,z^{'})=\sum\limits_{i,j=1}^nH_{ij}(z_0)z_iz_j,
\ee
where $H_{ij}=H_{ji}$ is a complex-valued matrix such that $\Im H$ is positive definite. It follows from (\ref{phase0}) and (\ref{var2}) that
\be\label{var22}
\sum\limits_{k,l=0}^n(\pl_{ij}^2\ol g^{kl}\pl_k\var\pl_l\var+2\ol g^{kl}\pl^3_{kij}\var\pl_l\var+2\ol g^{kl}\pl^2_{ki}\var\pl_{lj}^2\var+4\pl_i\ol g^{kl}\pl^2_{jk}\var\pl_l\var)|_{\beta}=0.
\ee
By the choices of $\var_0, \var_1$ and $\var_2,$ (\ref{var22}) implies that
$$(\pl^2_{ij}\ol g^{11}+2\ol g^{10}\pl^3_{0ij}\var+2\sum\limits_{k=2}^n\pl^2_{ki}\var\pl^2_{kj}\var)|_{\beta}=0.$$
We finally obtain the following Riccati equation for $H(z_0),$ namely,
\be\label{riccati}
\frac{d}{dz_0}H+HAH+B=0,\quad H(s)=H_0,\quad {\rm with}\quad \Im H_0>0\ {\rm and}\ z_0\in (a_0,b_0),
\ee
where $B=\frac14\pl_{ij}^2\ol g^{11},$ $s=\sqrt2t_-$ and the components of $A=(A_{ij})$ satisfy
$$\left\{ \begin{array}{l}
A_{11}=0,\\
A_{ii}=2,\quad i=2,\cdots,n,\\
A_{ij}=0,\quad {\rm otherwise}.
\end{array} \right.$$
For the above Riccati equation, we have
\begin{lem}\label{soriccati}\cite[Section 8]{AKYK}
The Riccati equation (\ref{riccati}) admits a unique solution. The solution $H$ is
symmetric and $\Im(H(z_0))>0$ for all $z_0\in (a_0,b_0).$ We have $H(z_0)=Z(z_0)Y^{-1}(z_0),$ where the matrix
valued functions $Z(z_0), Y(z_0)$ solve the first order linear system
$$\frac{dZ}{dz_0}=-BY,\quad \frac{dY}{dz_0}=AZ,\quad {\rm subject\ to}\quad Y(s)=I,\quad Z(s)=H_0.$$
Moreover, the matrix $Y(z_0)$ is non-degenerate on $(a_0,b_0)$, and there holds
$${\rm det}(\Im H(z_0))\cdot |{\rm det}(Y(z_0))|^2={\rm det}(\Im H_0).$$
\end{lem}

For the case where $|\Theta|=3,4\cdots,$ the polynomials $\var_j$ of higher degree are constructed analogously. We omit the details.
\vspace{0.1cm}

\noindent {\bf Construction of the amplitude.}\quad Let us consider the transport equation (\ref{transport0}). Let $|\Theta|=0.$ It follows from (\ref{nullg}) that
$b\var_t=\frac{-b}{\sqrt 2}$ and
$$\Del_{\ol g}\var=\sum\limits_{k,l=0}^n\ol g^{kl}D_{\ol g}^2\var(\pl_k,\pl_l)=\sum\limits_{k,l=0}^n\ol g^{kl}\pl_{kl}^2\var={\rm Tr}(AH)\quad {\rm on}\ \beta.$$
Therefore, the transport equation (\ref{transport0}) reduces to
\be\label{eqntrans}
2\pl_{z_0}a_{0,0}+[{\rm Tr}(AH)-\frac{b}{\sqrt 2}]a_{0,0}=0,\quad z_0\in (a_0,b_0).
\ee
Notice that
$${\rm Tr}(AH)={\rm Tr}(A(z_0)Z(z_0)Y^{-1}(z_0))={\rm Tr}(\frac{dY}{dz_0}Y^{-1}(z_0))=\frac{d}{dz_0}\log(\det Y(z_0)),$$
then
\be\label{solutrans}
a_{0,0}(z_0)=(\det(Y(z_0)))^{-\frac12}e^{\frac{1}{2\sqrt2}\int_s^{z_0}b(\tau,0) d\tau},\quad z_0\in (a_0,b_0)
\ee
is a solution to (\ref{eqntrans}). The subsequent terms $a_{k,0}$ can be constructed by solving some linear ODEs of first order. We refer to \cite{AFLO} for more details.
\vspace{0.2cm}

\noindent {\bf Construction of the remainder terms.}\quad By a similar proof with \cite[Lemma 2]{AFLO}, the Gaussian beam has the following property.
\begin{lem}\label{asyrh}
Let $u_\sigma=e^{\i\sigma\var}a$ be an approximate Gaussian beam of order $N$ along the null geodesic $\beta.$ Then for all $\sigma>0$
$$||\mathcal L_{b,q}u_\sigma||_{H^{k}(\mathcal M)}\le C\sigma^{-K},$$
where $K=\frac{N+1-k}{2}-1.$
\end{lem}
Based on the above lemma and the Sobolev embedding theorem, for sufficient large $N,$ the remainder term $r_{\sigma}$ satisfies the estimate (cf. \cite{Zhaijian}, \cite[Proposition 2.2]{LOPG})
\be\label{asyrsig}
||r_\sigma||_{H^{k+1}(\mathcal M)}\le C\sigma^{-K}\Rightarrow ||r_\sigma||_{C(\mathcal M)}\le C\sigma^{-\frac{n+1}{2}-2}.
\ee
\subsection{Proofs of Theorem 1.4}
Let us introduce some basic notations on the geodesic ray transform, which are explicitly discussed in, e.g., \cite{AFYK,AJYL,PG1}.
The following contents are mainly from \cite{AFYK}, we present here for completeness.

Let $S\Om\in T\Om$ be the unit sphere bundle of $(\Om,g)$, and by $\gamma(\cdot ;x,v)$ the geodesic with initial data $(x,v)\in S\Om.$ For all $(x,v)\in S\Om^{{\rm int}},$ we define the exist times as
$$\tau_{\pm}(x,v)=\inf\{r>0: \gamma(\pm r; x,v)\in \Ga\}.$$
Assume  that $(\Om,g)$ is simple, then $\tau_{\pm}<{\rm Diam}(\Om).$
Define
$$\pl_{\pm}S\Om=\{(x,v)\in S\Om: x\in \Ga,\quad \pm \<v,\nu(x)\>_g>0\}.$$
All geodesics in $\Om^{{\rm int}}$ can be parametrized by $\gamma(\cdot; x,v)\subset \Om$ for $(x,v)\in \pl_-S\Om.$ The geodesic ray transform on $(\Om,g)$ is defined for $f\in C^{\infty}(\Om)$ by
$$\mathcal{I}f(x,v)=\int_0^{\tau_+(x,v)}f(\gamma(r;x,v))dr\quad {\rm for}\ (x,v)\in\pl_-S\Om.$$
Let $\beta$ be a null geodesic (also called light ray).
By the product structure of the Lorentzian manifold $\R\times \Om$, we can parametrize the null geodesic $\beta$ as
$$\beta(r;s,x,v)=(r+s,\gamma(r;x,v)),\quad \forall (s,x,v)\in \R\times\pl_-S\Om.$$
Then, all the null geodesics $\beta$ through $\beta(\cdot ;s,x,v)$ with $(s,x,v)\in \R\times \pl_-S\Om$ over their maximal intervals $[0,\tau_+(x,v)]$ can be identified.
The so called light ray transform on $\R\times \Om$ can be defined as
$$\mathcal{G}f(s,x,v)=\int_0^{\tau_+(x,v)}f(r+s,\gamma(r;x,v))dr\quad \forall (s,x,v)\in \R\times\pl_-S\Om.$$


By \cite[Proposition 1.3]{AFYK} and the discussions in \cite[Section 2.2]{AJYL}, together with \cite[Theorem 1.6]{GMGH}, we have
\begin{pro}\label{lightray}
Suppose that either $(\Om,g)$ is simple or $(\Om,g)$ satisfies the foliation condition. Let $f\in C^1(\mathcal M))$ vanish on the set $\mathcal{M}\backslash\mathcal{E}.$ Then
$\mathcal{G}f=0$ for all maximal null geodesic $\beta\subset \mathcal D$ implies $f=0$.
\end{pro}
\begin{rem}
{\rm  For the non-simple case, such ray transform has been also  discussed e.g., in \cite{PG2} and references therein.}
\end{rem}

We begin with the first order linearized equation
\be\label{systemu1}
\left\{ \begin{array}{l}
\mathcal L_{b_j,q_j}u_j=u_{jtt}-\Del_gu_j+b_ju_{jt}+q_ju_j=0,\quad (t,x)\in \mathcal M,\\
u_j(t,x)=h_1(t,x),\qquad (t,x)\in \Sigma,\\
u_j(0,x)=u_{jt}(0,x)=0,\quad x\in \Om,
\end{array} \right.
\ee
where $q_j=f_{ju}(x,t,0).$
By the definition of $\mathcal L_{b,q},$ we have
$$\mathcal{L}_{b,q}^*u=u_{tt}-\Del_gu-bu_t+(q+b_t)u,$$
where $\mathcal L_{b,q}^*$ denotes the formal adjoint of $\mathcal L_{b,q}$ with respect to the $L^2(\mathcal M)$ inner product.
The formal adjoint system of (\ref{systemu1}) with $j=1$ is given by
\be\label{systemjoint}
\left\{ \begin{array}{l}
\mathcal L_{b_1,q_1}^*v=v_{tt}-\Del_g v-b_1v_t+(q_1+b_{1t})v=0,\ (t,x)\in \mathcal M,\\
v(T,x)=v_t(T,x)=0,\qquad x\in \Om,
\end{array} \right.
\ee

Based on the above constructions of Gaussian beams, we seek such solutions for systems (\ref{systemu1}) and (\ref{systemjoint}), respectively. More precisely, we let
$$u_2=e^{\i\sigma\var}a_1+r_{1\sigma}=e^{\i\sigma\var}\sigma^{\frac n4}\chi(\frac{|z^{'}|}{\delta^{'}})a_{10,0}+r_{1\sigma},$$ $$v=e^{-\i\sigma\ol\var}\ol a_2+r_{2\sigma}=e^{\i \sigma\var}\sigma^{\frac n4}\chi(\frac{|z^{'}|}{\delta^{'}})a_{20,0}+r_{2\sigma},$$
where $\ol\cdot$ means the conjugate of $\cdot.$
Then $r_{1\sigma}$ and $r_{2\sigma}$ respectively solve
\be\label{systemr1h}
\left\{\begin{array}{l}
\mathcal L_{b_2,q_2}r_{1\sigma}=-\mathcal L_{b_2,q_2}(e^{\i\sigma\var}a_1)\quad {\rm in}\ \mathcal M,\\
r_{1\sigma}=0 \qquad {\rm on}\ \Sigma,\\
r_{1\sigma}(0)=r_{1\sigma t}(0)=0 \quad {\rm in}\ \Om,
\end{array}\right.
\ee
and
\be\label{systemr2h}
\left\{\begin{array}{l}
\mathcal L^*_{b_1,q_1}r_{2\sigma}=-\mathcal L^*_{b_1,q_1}(e^{-\i\sigma\ol\var}\ol a_2)\quad {\rm in}\ \mathcal M,\\
r_{2\sigma}=0 \qquad {\rm on}\ \Sigma,\\
r_{2\sigma}(T)=r_{2\sigma t}(T)=0 \quad {\rm in}\ \Om,
\end{array}\right.
\ee
According to (\ref{lbq}), we have
$$\mathcal L_{b_2,q_2}(e^{\i\sigma\var}a_1)=e^{\i\sigma\var}[\mathcal L_{b_2,q_2}a_1+\sigma^2(\mathcal S\var)a_1+\i\sigma\mathcal T_{b_2}(a_1,\var)],$$
and
$$\mathcal L^*_{b_1,q_1}(e^{-\i\sigma\ol\var}\ol a_2)=e^{-\i\sigma\ol\var}[\mathcal L_{b_1,q_1}^*\ol a_2+\sigma^2(\mathcal S\ol\var)\ol a_2-\i\sigma\ol{\mathcal T_{-b_1}( a_2,\var)}].$$
As we have discussed in Section 4.1, we choose
$$a_{10,0}(z_0)=(\det Y(z_0))^{-\frac12}e^{\frac{1}{2\sqrt2}\int_s^{z_0}b_2(\tau,0) d\tau},\quad z_0\in (a_0,b_0),$$ and
$$a_{20,0}(z_0)=(\det Y(z_0))^{-\frac12}e^{-\frac{1}{2\sqrt2}\int_s^{z_0}b_1(\tau,0) d\tau},\quad z_0\in (a_0,b_0).$$

Clearly, $\mathcal L^*_{b_1,q_1}(e^{-\i \sigma\ol\var}\ol a_2)$ and $\mathcal L^*_{b_1,q_1}(e^{-\i \sigma\ol\var}\ol a_2)$ are compactly supported in a small tubular region around the null geodesic where the Fermi coordinates are well defined.

The following lemma shows that the remainder terms $r_{1\sigma}$ and $r_{2\sigma}$ vanish as $\sigma\rightarrow +\infty.$
\begin{lem}\label{remainder0}\cite{AJYL}
Let the remainder terms $r_{1\sigma}$ and $r_{2\sigma}$ respectively solve (\ref{systemr1h}) and (\ref{systemr2h}). Then
$r_{j\sigma}\in C([0,T];H_0^1(\Om))\cap C^1([0,T];L^2(\Om)),$ and
$$\mathop{\lim}\limits_{\sigma\to\infty}(||r_{j\sigma}||_{L^2(\mathcal M)}+\sigma^{-1}||r_{j\sigma}||_{H^1(\mathcal M)})=0,\quad {\rm for}\quad j=1,2.$$
\end{lem}
We are now in a position to prove Theorem 1.4.

We only prove Theorem 1.4 for the case where $(M,g)$ is simple. If $(\Om,g)$ satisfies the foliation condition, as in \cite{GUAV} (see also \cite{Zhaijian}), for any point $q\in \Ga$,
there exists a wedge-shaped neighborhood $O_q\subset\Om$ of $q$ such that any geodesic in $(O_q, g)$ has no conjugate points. Therefore, we can now recover $f_z^{k}(x,t,0)$ for $k\ge 3.$ Then the foliation
condition allows a layer stripping scheme to recover the coefficients in the whole domain. However, the recovery of $f_z(x,t,0)$ and $f_z^{(2)}(x,t,0)$ is quite different, which needs the inversion of some new ray transforms.
\\

\noindent{\bf Proof of Theorem 1.4.}\quad We divide the proof into three steps.

{\bf Step 1.}
Let $h_1=e^{\i\sigma\var}a_1|_{\Sigma}$ in (\ref{systemu1}). Let
$$w=u_1-u_2,\quad b=b_1-b_2,\quad q=q_1-q_2.$$ Then
\be\label{systemubq}
\left\{ \begin{array}{l}
w_{tt}-\Del_gw+b_1w_t+q_1w=-(bu_{2t}+qu_2),\quad (t,x)\in \mathcal M,\\
w(t,x)=0,\qquad (t,x)\in \Sigma,\\
w(0,x)=w_t(0,x)=0,\quad x\in \Om.
\end{array} \right.
\ee
Notice that $v$ solves system (\ref{systemjoint}) and $\Lambda_{b_1,q_1}(h_1)=\Lambda_{b_2,q_2}(h_1)$ on $\Sigma$ implies $\pl_{\nu}w|_{\Sigma}=0,$ we multiply the first equation of (\ref{systemubq}) by $v$ and integrate over $\mathcal M$ to obtain
\be
(\mathcal L_{b_1,q_1}w,v)_{L^2(\mathcal M)}=(w,\mathcal L^*_{b_1,q_1}v)_{L^2(\mathcal M)}=0\Rightarrow \int_{\mathcal M}(bu_{2t}+qu_2)vdV_{\ol g}=0.
\ee
Here $dV_{\ol g}=|\ol g|^{\frac12}dt\wedge dx$ denotes the volume form of the metric $\ol g=-dt^2+g.$
Recalling that $u_2=e^{\i \sigma\var}a_1+r_{1\sigma}$ and $v=e^{-\i \sigma\ol\var}\ol a_2+r_{2\sigma}$ in a neighborhood of the null geodesic $\beta,$ we have
\beq\label{interbq}
0&&=\i \sigma\int_{\mathcal M} b\var_ta_1\ol a_2e^{-2\sigma{\rm Im}\var} dV_{\ol g}+\int_{\mathcal M}qa_1\ol a_2e^{-2\sigma{\rm Im}\var}dV_{\ol g}\nonumber\\
&&\quad+\int_{\mathcal M} b[a_{1t}\ol a_2e^{-2\sigma{\rm Im}\var}+(\i \sigma\var_ta_1+a_{1t})r_{2\sigma}e^{\i \sigma\var}+e^{-\i \sigma\ol\var}\ol a_2r_{1\sigma t}+r_{1ht}r_{2\sigma}]dV_{\ol g}\nonumber\\
&&\quad+\int_{\mathcal M} q(e^{\i \sigma\var}a_1r_{2\sigma}+e^{-\i\sigma\ol\var}\ol a_2r_{1\sigma}+r_{1\sigma}r_{2\sigma})dV_{\ol g}
\eeq
It then follows from Lemma \ref{remainder0} that
\be\label{limitb}
\mathop{\lim}\limits_{\sigma\to\infty}\int_{\mathcal M}b\var_ta_1\ol a_2e^{-2\sigma{\rm Im}\var}dV_{\ol g}=0.
\ee
Since the functions $a_1,a_2$ are supported in a small tubular neighborhood of the null geodesic $\beta$, the integrand in (\ref{limitb}) is supported near $\beta.$ Therefore we can use the Fermi coordinates $z=(z_0,z_1,z^{''})$ to compute the limit. Recall that
$$a_1\ol a_2=\sigma^{\frac n2}a_{10,0}\ol a_{20,0}\chi^2(\frac{|z^{'}|}{\delta^{'}})+\mathcal O(\sigma^{-1})=\sigma^{\frac n2}|\det Y(z_0)|^{-1}e^{-\frac{1}{2\sqrt2}\int_s^{z_0}b(\tau,0)d\tau}+\mathcal O(\sigma^{-1}).$$
Moreover, by \cite[Lemma 3.6]{AJYL}, we know that $\var_t$ does not vanish in $\ol{\mathcal V}.$ Notice that $b,q=0$ on the set $\mathcal M_1\backslash \mathcal M.$ Then (\ref{limitb}) yields
\be\label{geobq}
\mathop{\lim}\limits_{\sigma\to\infty}\sigma^{\frac n2}\int_{a_0}^{b_0}\int_{|z^{'}|<\delta^{'}}\chi^2(\frac{|z^{'}|}{\delta^{'}})e^{-2\sigma{\rm Im}\var}|\det Y(z_0)|^{-1}b(z_0,z^{'})e^{-\frac{1}{2\sqrt2}\int_s^{z_0}b(\tau,0)d\tau}dz_0\wedge z^{'}=0.\nonumber
\ee

We proceed to calculate the term
$$\sigma^{\frac n2}\int_{|z^{'}|<\delta^{'}}\chi^2({\frac{|z^{'}|}{\delta^{'}}})b(z_0,z^{'})e^{-2\sigma{\rm Im}\var}dz^{'}=\int_{\R^n}e^{-2\sigma x^TPx}\eta(x)dx,$$
where $\eta(x)=\chi^2(\frac{|x|}{\delta^{'}})b(z_0,x)$ is a smooth function with compact support $B(0,\delta^{'})$, $P={\rm Im}H(z_0)$ is a positive-definite matrix.
By the following well-known formula
$$\mathcal F(e^{-\frac12x^TPx})(\xi)=\frac{(2\pi)^{\frac n2}}{(\det P)^{\frac12}}e^{-\frac12\xi^TP^{-1}\xi},$$
where $\mathcal F$ denotes the Fourier transform, we have
$$\mathcal F(e^{-2\sigma x^TPx})(\xi)=\frac{(2\pi)^{\frac n2}}{(\det P)^{\frac12}(4\sigma)^{\frac n2}}e^{-\frac{1}{8\sigma}\xi^TP^{-1}\xi}.$$
Equivalently, we have
$$\mathcal F[e^{-\frac12x^T(\frac{1}{4\sigma}P^{-1})x}]=\frac{(2\pi)^{\frac n2}}{\det(\frac{1}{4\sigma}P^{-1})^{\frac12}}e^{-2\sigma x^TPx}.$$
Since
$$\int_{\R^n}\mathcal F[f](x)g(x)dx=\int_{\R^n}f(x)\mathcal F[g](x)dx,\quad {\rm for}\quad f,g\in L^p(\R^n),\quad p\in (1,+\infty),$$
we have
\beq
\int_{\R^n}e^{-\frac{1}{8\sigma}\xi^TP^{-1}\xi}\mathcal F\eta(\xi)d\xi&&=\int_{\R^n}\mathcal F[e^{-\frac{1}{8\sigma}x^TP^{-1}x}]\eta(x)dx\nonumber\\
&&=(2\pi)^{\frac n2}(4\sigma)^{\frac n2}(\det P)^{\frac12}\int_{\R^n}e^{-2\sigma x^TPx}\eta(x)dx.\nonumber
\eeq
Let $\PP(\xi): \R^n\to \mathbb C$ be a measurable function. The well known theory of pseudo-differential operators tells us
$$\mathcal F[\PP(D)u](\xi)=\PP(\xi)\hat u(\xi)=\PP(\xi)\mathcal F[a](\xi),\quad \PP(\xi)\in S^m,$$
where $D$ is the differential operator, $S^m$ is the symbol class of order $m,$ and
$$\PP(D)u=\frac{1}{(2\pi)^n}\int_{\R^n}e^{\i x\cdot\xi}\PP(\xi)\hat u(\xi)d\xi.$$
Therefore
\beq
&&\int_{\R^n}e^{-2hx^TPx}\eta(x)dx=\frac{1}{(2\pi)^{\frac n2}(4h)^{\frac n2}(\det P)^{\frac12}}\int_{\R^n}e^{-\frac{1}{8h}\xi^TP^{-1}\xi}\mathcal F \eta(\xi)d\xi\nonumber\\
&&\quad=\frac{1}{(2\pi)^{\frac n2}(4\sigma)^{\frac n2}(\det P)^{\frac12}}\sum\limits_{k=0}^{+\infty}\frac{1}{k!}(-\frac{1}{8\sigma})^k\int_{\R^n}(\xi^TP^{-1}\xi)^k\mathcal F \eta(\xi)d\xi\nonumber\\
&&\quad=\frac{1}{(4\sigma)^{\frac n2}(\det P)^{\frac12}}\sum\limits_{k=0}^{+\infty}\frac{1}{k!}(-\frac{1}{8\sigma})^k({\mathcal{P}_{P^{-1}}(D))^k\eta(0)}\nonumber\\
&&\quad=\frac{1}{(4\sigma)^{\frac n2}(\det P)^{\frac12}}(\eta(0)+\mathcal O (\sigma^{-1})),\nonumber
\eeq
where
${\mathcal{P}_{P^{-1}}(D)\eta(\xi)}$ is defined by
$$(\xi^TP^{-1}\xi)^k\mathcal F\eta(\xi)=\mathcal F[({\mathcal{P}_{P^{-1}}(D))^k\eta(\xi)}].$$
By Lemma \ref{soriccati}, we have
$$|\det Y(z_0)|^{-1}=(\det\Im H(z_0))^{\frac12}(\det\Im H_0)^{-\frac12}$$
Then
\beq\label{detbb}
&&\mathop{\lim}\limits_{\sigma\to\infty}\sigma^{\frac n2}\int_{|z^{'}|<\delta^{'}}\chi^2({\frac{|z^{'}|}{\delta^{'}}})e^{-2\sigma{\rm Im}\var}|\det Y(z_0)|^{-1}e^{-\frac{1}{2\sqrt2}\int_s^{z_0}b(\tau,0)d\tau}dz^{'}\nonumber\\
&&\qquad=\frac{1}{4^n(\det\Im H_0)^{-\frac12}}b(z_0,0)e^{-\frac{1}{2\sqrt2}\int_s^{z_0}b(\tau,0)d\tau}.
\eeq
Combine (\ref{limitb}) with (\ref{geobq}), we find that
\be
-2\sqrt2\int_{a_0}^{b_0}\frac{d}{dz_0}e^{-\frac{1}{2\sqrt2}\int_s^{z_0}b(\tau,0)d\tau}dz_0=\int_{a_0}^{b_0}b(z_0,0)e^{-\frac{1}{2\sqrt2}\int_s^{z_0}b(\tau,0)d\tau}dz_0=0.
\ee
Therefore
$$\int_{a_0}^{b_0}b(z_0,0)dz_0=0,$$
which implies that
$$\int_{\beta}b(\beta)=0$$
holds for any maximal  null geodesic in $\R\times \Om$.
Then Proposition \ref{lightray} is applied to obtain $b=0,$ which implies that $b_1=b_2$ in $\mathcal E.$

We return to the equation (\ref{interbq}) with $b=0$ to have
$$\mathop{\lim}\limits_{\sigma\to\infty}\int_{\mathcal M}qa_1\ol a_2e^{-2\sigma{\rm Im}\var}dV_{\ol g}=0.$$
By a similar argument, this reduces
\beq
&&\mathop{\lim}\limits_{\sigma\to\infty}\sigma^{\frac n2}\int_{a_0}^{b_0}\int_{|z^{'}|<\delta^{'}}q\chi^2(\frac{|z^{'}|}{\delta^{'}})|\det Y(z_0)|^{-1}e^{-2\sigma{\rm Im}\var}dz_0\wedge dz^{'}=0.
\eeq
Then
$$\int_{\beta}q(\beta)=0$$
holds for all maximal null geodesics $\beta$ in $\R\times\Om.$ Together with Proposition \ref{lightray}, we conclude that $q=0.$
\\

{\bf Step 3.}\quad Since $f_z^{(2)}$ is known, we proceed with the  third and higher order linearization.  Let
$$W^{(123)}=\frac{\pl^3}{\pl\epsilon_1\pl\epsilon_2\pl\epsilon_3}\Big|_{\epsilon=0}u(f),\quad W^{(ij)}=\frac{\pl^2}{\pl\epsilon_i\pl\epsilon_j}\Big|_{\epsilon=0}u(f),\quad 1\le i,j\le 3.$$
Set $m(x,t)=f_{uuu}(x,t,0).$ Let $\Sigma(3)$ be the permutation group of $\{1,2,3\}.$ Then
\be
\left\{ \begin{array}{l}
\mathcal L_{b,q}W^{(123)}+\frac{f_{uu}(x,t,0)}{2}\sum\limits_{\zeta\in \Sigma(3)}W^{(\zeta(1)\zeta(2))}v_{\zeta(3)}+mv_1v_2v_3=0\quad {\rm in}\ \mathcal M,\\
W^{(123)}=0\qquad {\rm on}\ \Sigma,\\
W^{(123)}(x,0)=W^{(123)}_t(x,0)=0\quad {\rm in}\ \Om,
\end{array} \right.
\ee
where $$\mathcal L_{b,q}W^{(123)}=W^{(123)}_{tt}-\Del_gW^{(123)}+bW^{(123)}_t+qW^{(123)},$$
and for each $k=1,2,3,$ $v_k=\frac{\pl}{\pl\epsilon_k}\Big|_{\epsilon=0}u(f)$ solves
\be\label{sysuj}
\left\{ \begin{array}{l}
v_{ktt}-\Del_gv_k+bv_{kt}+qv_k=0\quad {\rm in}\ \mathcal M,\\
v_k=h_k\qquad {\rm on}\ \Sigma,\\
v_k(x,0)=v_{kt}(x,0)=0\quad {\rm in}\ \Om.
\end{array} \right.
\ee

Let $v_0$ solve the following adjoint system of (\ref{sysuj})
\be
\left\{ \begin{array}{l}
v_{0tt}-\Del_gv_0-bv_{0t}+(q+b_t)v_0=0\quad {\rm in}\ \mathcal M,\\
v_0=h_0\qquad {\rm on}\ \Sigma,\\
v_0(x,T)=v_{0t}(x,T)=0\quad {\rm in}\ \Om.
\end{array} \right.
\ee
Integrating by parts over $\mathcal M$ yields
\beq
&&\int_{\Sigma}v_0\frac{\pl^3}{\pl\epsilon_1\pl\epsilon_2\pl\epsilon_3}\Big|_{\epsilon=0}\Lambda(\epsilon h_1+\epsilon_2h_2+\epsilon_3h_2)d\Sigma\nonumber\\
&&=\int_{\mathcal M} mv_1v_2v_3v_0+\int_{\mathcal M} \frac{f_{uu}(x,t,0)}{2}\sum\limits_{\zeta\in \Sigma(3)}W^{(\zeta(1)\zeta(2))}v_{\zeta(3)}v_0dgdt.
\eeq
Since $f_{2uu}(x,t,0)$ is known. Then the Dirichlet to Neumann map $\Lambda$ determines
\be\label{int0123}
\int_{\mathcal M}mv_1v_2v_3v_0dV_{\ol g}.
\ee
We will use special solutions $v_1,v_2,v_3,v_0$ in the above identity. Concretely, we shall use the following Gaussian beam solutions $e^{\i\sigma\var}a+r_{\sigma} $
constructed in section 4.1.

For given $p=(t_0,x_0)\in \mathcal E,$ we choose local coordinates such that $\ol g$ coincides with the standard Minkowski metric at $p.$ Define the light cone at $p$ as $$\mathcal C(p)=\{(t,X)\in T_p\mathcal M: t^2=|X|_g^2\}.$$
Similar to \cite[Lemma 1]{XMLG} (see also, \cite[section 3.2]{Zhaijian}, \cite[Section 5]{AFLO}), we can assume without loss of generality that
$\zeta_0,\zeta_1\in \mathcal C(p)$ satisfying
$$\zeta_0=(1,-\sqrt{1-\theta^2},\theta,0,\cdots,0),\quad \zeta_1=(1,1,0,\cdots,0)$$ for some $\theta\in [0,1].$ Taking $\t\theta>0$ small and introduce
$$\zeta_2=(1,\sqrt{1-\t\theta^2},\t\theta,0,\cdots,0),\quad \zeta_3=(1,\sqrt{1-\t\theta^2},-\t\theta,0,\cdots,0).$$
By \cite[Lemma 1]{XMLG}, $\zeta_0,\zeta_1$ are linear-independent, and there are constants $k_0,k_1,k_2,k_3$ such that
$$k_0\zeta_0+k_1\zeta_1+k_2\zeta_2+k_3\zeta_3=0.$$
Denote $\beta_k$ to be the null geodesic with cotangent vector $\zeta_k$ and $p.$
Taking $$v_k=e^{\i \sigma k_k\var_k}a_k+r_{k\sigma}\quad {\rm for}\quad k=0,1,2,3$$ as the Gaussian beams concentrating near the null geodesic $\beta_k.$
Notice that the manifold $(\Om,g)$ is simple, the null geodesic $\beta_k$ $(k=0,1,2,3)$ can intersect only at $p.$

Inserting $v_k=e^{\i \sigma k_k\var_k}a_k+r_{k\sigma}$ into (\ref{int0123}), with estimate (\ref{asyrsig}), the Dirichlet-to-Neumann map determines
\beq\label{asyphase}
&&\sigma^{\frac{n+1}{2}}\int_{\mathcal M} mv_0v_1v_2v_3dV_{\ol g}\nonumber\\
&&=\sigma^{\frac{n+1}{2}}\int_{\mathcal M}me^{\i\sigma(k_0\var_0+k_1\var_1+k_2\var_2+k_3\var_3)}a_0a_1a_2a_3dV_{\ol g}+\mathcal O(\sigma^{-1}).
\eeq
Clearly, the product $a_0a_1a_2a_3$ is supported in a neighborhood of $p.$ We introduce a lemma to deal with the above integral.
\begin{lem}\label{vectors}\cite[Lemma 5]{AFLO}. The function
$$S:=k_0\var_0+k_1\var_1+k_2\var_2+k_3\var_3$$
is well-defined in a neighborhood of $p$ and

(1) $S(p)=0;$

(2) $D^{\ol g}S(p)=0;$

(3) $\Im S(p_1)\ge cd(p_1,p)$ for $p_1$ in a neighborhood of $p,$ where $c>0$ is a constant.

\end{lem}

Based on the above lemma, applying the stationary phase (e.g., see \cite[Theorem 7.75]{Homander}) to (\ref{asyphase}), we
$$c\sigma^{\frac{n+1}{2}}\int_{\mathcal M} mv_0v_1v_2v_3dV_{\ol g}=m(p)(a_0a_1a_2a_3)(p)+\mathcal O(\sigma^{-1}),$$
where $c$ denotes some explicit constant.
Thus, the Dirichlet-to-Neumann map determines $m(p).$

For the recovery the higher order coefficients $f^{(k)}_u(x,t,0)$ for $k\ge 4,$ we can achieve this by induction. We
refer to \cite[Section 4]{Zhaijian} for such an operation and omit the details. Therefore, the proof of Theorem 1.4 is complete. \qedsymbol

\begin{rem}
{\rm We notice that the recovery of $q$ and $b$ is much different from that of higher order terms $f^{(k)}_u$ ($k\ge 3$). Lemma \ref{remainder0} is a special case of Lemma \ref{asyrh}, and the Gaussian beam of order $N=0$ is enough for the linear inverse problem of recovering $b$ and $q$.

It seems that we can not recover $f_{uu}$ by the same method as that for terms $f_u^{(k)}$ for $k=1\ {\rm or}\ k\ge 3$. One of the reason is that we can not choose three time-like vectors $\zeta_0,\zeta_1,\zeta_2$ such that $\zeta_0, \zeta_1$ are linear-dependent but $\zeta_0,\zeta_1,\zeta_2$ are linear-dependent. We mention that, in \cite{ZU}, due to the presence of two different matrics $g_P$ and $g_S,$ the authors have chosen three different vectors satisfying the above Lemma \ref{vectors}.  Therefore, they proved the unique recovery of coefficients by the Gaussian beams in place of the distorted plane waves method used in e.g., \cite{Zhaijian}.}
\end{rem}

\section*{Appendix \ Well-posedness}
We prove the well-posendess result to (1.1).
We begin with the following linear wave equation
\be\label{systemb}
\left\{ \begin{array}{l}
u_{tt}-\Del_gu=F\quad {\rm in}\ \mathcal M,\\
u=h\qquad {\rm on}\ \Sigma,\\
u(x,0)=u_0, u_t(x,0)=u_1(x)\quad {\rm in}\ \Om.
\end{array} \right.
\ee
Let $$(u_0,u_1)\in H^{m+1}(\Om)\times H^m(\Om),\quad h\in H^{m+1}(\Sigma),\quad F\in L^1([0,T];H^m(\Om))$$ with $\pl_t^kF(x,t)\in L^1([0,T];H^{m-k}(\Om))$ for $k=0, 1,\cdots m.$
Moreover, we assume that the compatibility conditions hold up to order $m,$ which are given by
\be\label{comcondition}
\left\{ \begin{array}{l}
h(x,0)=u_0(x)|_{\Ga},\quad h_t(x,0)=u_1(x)|_{\Ga},\quad h_{tt}(x,0)=[\Del_gu_0(x)+F(x,0)]|_{\Ga},\vspace{1.0ex}\\
\pl_t^kh(x,0)=\pl_t^ku(x,0)|_{\Ga},\ k=3,\cdots, m.
\end{array} \right.
\ee
According to \cite[Theorem 2.45]{AKYK}, for system {\ref{systemb}} with the above conditions, we have
\begin{lem}\label{welllin}
Let $m$ be a positive integer and $T>0.$ Then system (\ref{systemb}) admits a unique solution
$u\in E^{m+1}$ and $\pl_{\nu}u\in H^m(\Sigma).$ Moreover, the dependence of $u$ and $\pl_{\nu}u$ on $u_0,u_1,h,F$ is continuous in the corresponding spaces, i.e.,
\beq
&&||u||_{E^{m+1}}+||\pl_{\nu}v||_{H^m(\Sigma)}\nonumber\\
&&\le C(T)(||u_0||_{H^{m+1}(\Om)}+||u_1||_{H^m(\Om)}+||h||_{H^{m+1}(\Sigma)}+||F||_{X_m}),
\eeq
where $||F||^2_{X_m}=\sum\limits_{k=0}^m||\pl^k_tF||^2_{L^1([0,T];H^{m-k})}.$
\end{lem}
Assume further that $F\in E^m$ and $b,q\in C^m(\mathcal M)$ with $m>\frac n2.$ Assume that the compatibility conditions hold for (\ref{sysbbm}) up to order $m$. Applying lemma \ref{welllin}, and by the fact that $E^m$ is a Banach algebra when $m>\frac n2,$ we know that
\be\label{sysbbm}
\left\{ \begin{array}{l}
v_{tt}-\Del_gv+bv_t+qv=F\quad {\rm in}\ M,\\
v=h\qquad {\rm on}\ \Sigma,\\
v(x,0)=v_0, v_t(x,0)=v_1(x)\quad {\rm in}\ \Om.
\end{array} \right.
\ee
admits a unique solution $v\in E^{m+1}$ with $\pl_{\nu}v\in H^{m}(\Sigma).$ Moreover, the following energy estimate holds
\beq
&&||v||_{E^{m+1}}+||\pl_{\nu}v||_{H^m(\Sigma)}\nonumber\\
&&\le C(T)(||v_0||_{H^{m+1}(\Om)}+||v_1||_{H^m(\Om)}+||h||_{H^{m+1}(\Sigma)}+||F||_{X_m}).
\eeq

We now in a position to prove the  well-posedness of the nonlinear system (\ref{systembf}) with small initial data $(u_0,u_1)$ and small boundary data $h.$

Let $v$ solve the following non-homogeneous linear wave equation
\be\label{syshomo}
\left\{ \begin{array}{l}
v_{tt}-\Del_gv+bv_t+f_z(x,t,0)v=0\quad {\rm in}\ M,\\
v=h\qquad {\rm on}\ \Sigma,\\
v(x,0)=u_0, v_t(x,0)=u_1(x)\quad {\rm in}\ \Om.
\end{array} \right.
\ee
For $L>0$ given, let $$B_{m+1}(L)=\{(y_1,y_2)\in H^{m+1}(\Om)\times H^m(\Om): ||y_1||_{H^{m+1}}+||y_2||_{H^m}\le L\},$$
$$\mathcal N_{m+1}(L)=\{h\in H^{m+1}(\Sigma): ||h||_{H^{m+1}(\Sigma)}\le L\}\subset H^{m+1}(\Sigma).$$
Let $u_0,u_1\in B_{m+1}(\varepsilon_0/3),$ and $h\in \mathcal N_{m+1}(\varepsilon_0/3)$ for $\varepsilon_0$ sufficiently small. Then
\be\label{estv}
||v||_{E^{m+1}}+||\pl_{\nu}v||_{H^m(\Sigma)}\le C(T)\varepsilon_0.
\ee
For given $\delta_0>0$ sufficiently small,
let $$Z_{m+1}(\delta_0)=\{\hat w\in E^{m+1}: ||\hat w||_{E^{m+1}}\le \delta_0 \}\subset E^{m+1}.$$
Let $w$ solve the following homogeneous equation
\be\label{sysnon}
\left\{ \begin{array}{l}
w_{tt}-\Del_gw+bw_t+f_z(x,t,0)w+\sum\limits_{k=2}^{\infty}f^{(k)}(x,t,0)\frac{(v+\hat w)^k}{k!}=0\quad {\rm in}\ \mathcal M,\\
w=0\qquad {\rm on}\ \Sigma,\\
w(x,0)=w_t(x,0)=0\quad {\rm in}\ \Om,
\end{array} \right.
\ee
where $v\in E^{m+1}$ is the solution of (\ref{syshomo}) with estimate (\ref{estv}).
We define a map $$\mathcal A: Z_{m+1}(\delta_0)\rightarrow E^{m+1},$$ which sends the given $\hat w\in Z_{m+1}(\delta_0)$ to the solution of (\ref{sysnon}).
For any positive integers $k$ and $R,$ one has
$$||f^{(k)}(x,t,0)||_{E^{m+1}}\le \frac{k!}{R^k}\mathop{\sup}\limits_{|z|=R}||f(x,t,z)||_{E^{m+1}}.$$
By the a priori estimate of (\ref{sysnon}) and the property of Banach algebra, we have
\beq
&&||\mathcal A\hat w||_{E^{m+1}}\le C(T)\sum\limits_{k=2}^{\infty}\frac{1}{k!}||f^{(k)}(x,t,0)||_{E^{m+1}}||v+\hat w||^k_{E^{m+1}}\nonumber\\
&&\le C(T)\sum\limits_{k=2}^{\infty}\frac{1}{R^k}\mathop{\sup}\limits_{|z|=R}||f(x,t,z)||_{E^{m+1}}2^{k-1}(||v||^k_{E^{m+1}}+||\hat w||^k_{E^{m+1}})\nonumber\\
&&\le \frac{C(T)}{R}\mathop{\sup}\limits_{|z|=R}||f(x,t,z)||_{E^{m+1}}\sum\limits_{k=1}^{\infty}\frac{2^k}{R^k}(\delta_0\delta_0^k+\varepsilon_0\varepsilon_0^k)\nonumber\\
&&=\frac{2C(T)}{R}\mathop{\sup}\limits_{|z|=R}||f(x,t,z)||_{E^{m+1}}(\delta_0\frac{\delta_0}{R-2\delta_0}+\varepsilon_0\frac{\varepsilon_0}{R-2\varepsilon_0}).
\eeq
Let $R=2$ and $\varepsilon_0=\delta_0\le \frac12$ small enough such that
$$C(T)\mathop{\sup}\limits_{|z|=R}||f(x,t,z)||_{E^{m+1}}\delta_0\le \frac12.$$
Then we have
$$||\mathcal A\hat w||_{E^{m+1}}\le \delta_0,$$
which implies that $\mathcal A: Z_{m+1}(\delta_0)\rightarrow Z_{m+1}(\delta_0)$ is well defined.
Let
$$F(x,t,\hat w)=\sum\limits_{k=2}^{\infty}f^{(k)}(x,t,0)\frac{(v+\hat w)^k}{k!}.$$
Then
\beq
&&F(x,t,\hat w_1)-F(x,t,\hat w_2)\nonumber\\
&&=\sum\limits_{k=2}^{\infty}\frac{f^{(k)}(x,t,0)}{(k-1)!}\int_0^1[\tau(v+\hat w_1)+(1-\tau)(v+\hat w_2)]^{k-1}d\tau.\nonumber
\eeq
Therefore, taking $\delta_0$ small enough, we see that
\beq
&&||\mathcal A \hat w_1-\mathcal A\hat w_2||_{E^{m+1}}=||w_1- w_2||_{E^{m+1}}\nonumber\\
&&\le C(T)\sum\limits_{k=2}^{\infty}\frac{k}{R^k}\mathop{\sup}\limits_{|z|=R}||f(x,t,z)||_{E^{m+1}}(3\delta_0)^{k-1}||\hat w_1-\hat w_2||_{E^{m+1}}\nonumber\\
&&\le \frac12||\hat w_1-\hat w_2||_{E^{m+1}}.
\eeq
Thus, $\mathcal A: Z_{m+1}(\delta_0)\rightarrow Z_{m+1}(\delta_0)$ is a contraction. The Banach's fixed point theorem implies that $u=v+w\in E^{m+1}$ is a solution to system (\ref{systembf}). Moreover, we have
$$||u||_{E^{m+1}}+||u||_{C(\ol{\mathcal M})}+||\pl_{\nu}u||_{H^m(\Sigma)}\le C(||u_0||_{H^{m+1}(\Om)}+||u_1||_{H^m(\Om)}+||h||_{H^{m+1}(\Sigma)}),$$
where $C>0$ is independent of $u_0,u_1$ and $h.$

\end{document}